# On the Quaternionic Curves in the Semi-Euclidean Space $E_2^4$

Tülay SOYFİDAN, Mehmet Ali GÜNGÖR

**Abstract:** In this study, we try to semi-real quaternionic curves in the semi-Euclidean space $\mathbb{E}_2^4$. Firstly, we introduce algebraic properties of semi-real quaternions. And then, we give some characterizations of semi-real quaternionic involute-evolute curves in the semi-Euclidean space $\mathbb{E}_2^4$. Lastly, we illustrate some examples and draw their figures with Mathematica Programme.

## 1. Introductions

Quaternions, for the first time in 1843, was discovered by the Irish mathematician Sir William R. Hamilton. Hamilton wanted to generalize complex numbers to use geometric optics. Thus, quaternions, which are a more general form of complex numbers, was found by him, [1]. Quaternions, as well as having a normal vector algebra of finite rotations for calculated mathematical calculations of physical problems, they provides a simple and elegant especially for describing finite rotations in space. Also, they have many useful methods, such as Euler angles and in mechanics, for example, quaternionic formulation of equation of motion in the theory relativity.

As a set, the quaternions $\mathbb{Q}$ are coincide with $\mathbb{R}^4$, a four-dimensional vector space over the real numbers. The Serret-Frenet formulae for a quaternionic curves in $\mathbb{R}^3$ are introduced by K. Bharathi and M. Nagaraj. Moreover, they obtained the Serret-Frenet formulae for the quaternionic curves in $\mathbb{R}^4$ by the formulae in $\mathbb{R}^3$, [2]. Then, lots of studies have been published by using this studies. One of them is A. C. Çöken and A. Tuna's study [3] which they gave Serret-Frenet formulas, inclined curves, harmonic curvatures and some characterizations for a quaternionic curve in the semi-Euclidean spaces $\mathbb{E}_1^3$ and $\mathbb{E}_2^4$. Another is Gök et al.'s and Kahraman et al.s studies. They defined a new kind of slant helix, which they called $B_2$–slant helix and they gave some characterizations of this slant helix in $\mathbb{E}_2^4$, [4].

This study deals with special types of curves called involutes and evolutes. These are curves which can be associated to any given curve $\alpha$. An involute to $\alpha$ may be thought of as any curve which is always perpendicular to the tangent lines of $\alpha$. If you think of unwinding string tautly from around a curve, you will get a picture of an involute. An evolute of is a curve whose involute is $\alpha$, so evolutes and involutes are, in a sense, inverses to each other. An evolute of may also be thought of as the curve determined by the centers of curvature of $\alpha$. Many studies have been conducted on these curves. One of them Bükçü and Karacan's study which they generalized the involute and evolute curves of the spacelike curve $\alpha$ with a spacelike binormal in Minkowski 3-Space, [6].

The main purpose of this paper is to obtain some characterizations of semi-real quaternionic involute-evolute curves in semi-quaternionic space $\mathbb{Q}_v$. To do this, since it is a trivial task to write out the Serret-Frenet formulae of the curve in $\mathbb{E}_2^4$ using semi-real quaternions, firstly it is established some characterizations of semi-real quaternionic involute-evolute curves in $\mathbb{E}_2^4$. And some results for semi-real quaternionic $w-curves$ which has constant curvatures are given. Moreover, it is seen that the semi-real spatial quaternionic curves in $\mathbb{E}_1^3$ associated with semi-real quaternionic involute-evolute curves in $\mathbb{E}_2^4$ aren't semi-real quaternionic involute-evolute curves. Lastly, we illustrate some examples and draw their figures with Mathematica Programme.

## 2. Preliminaries

A semi-real quaternion is defined with $q = q_1 \boldsymbol{e}_1 + q_2 \boldsymbol{e}_2 + q_3 \boldsymbol{e}_3 + q_4$ (or $q = S_q + \boldsymbol{V}_q$ where the symbols $S_q = q_4$ and $\boldsymbol{V}_q = q_1 \boldsymbol{e}_1 + q_2 \boldsymbol{e}_2 + q_3 \boldsymbol{e}_3$ denote scalar and vector part of $q$, respectively) such that

$$\begin{aligned} i) \quad & \boldsymbol{e}_i \times \boldsymbol{e}_i = -\varepsilon_{e_i}, \qquad (1 \leq i \leq 3) \\ ii) \quad & \boldsymbol{e}_i \times \boldsymbol{e}_j = \varepsilon_{e_i} \varepsilon_{e_j} \boldsymbol{e}_k \quad \text{in } \mathbb{R}_1^3 \\ & \boldsymbol{e}_i \times \boldsymbol{e}_j = -\varepsilon_{e_i} \varepsilon_{e_j} \boldsymbol{e}_k \quad \text{in } \mathbb{R}_2^4 \end{aligned} \qquad (2.1)$$

where $(ijk)$ is an even permutation of $(123)$. For every $p, q \in \mathbb{Q}_v$, using these basic products we can now expand the product of two semi-real quaternion as

$$p \times q = S_p S_q + g(\boldsymbol{V}_p, \boldsymbol{V}_q) + S_p \boldsymbol{V}_q + S_q \boldsymbol{V}_p + \boldsymbol{V}_p \wedge_{\mathbb{L}} \boldsymbol{V}_q, \qquad (2.2)$$

where we have used the usual inner and cross products in semi-Euclidean space $\mathbb{R}_1^3$, [3]. A feature of semi-real quaternions is that the product of two semi-real quaternions is non-commutative. The conjugate of the semi-real quaternion $q$ is denoted by $\overline{q}$ and defined $\overline{q} = S_q - \boldsymbol{V}_q$ Thus, we define symmetric, non-degenerate valued bilinear form $h$ as follow:

$$\begin{aligned} h(p,q) &= \frac{1}{2} \left[ \varepsilon_p \varepsilon_q (p \times \overline{q}) + \varepsilon_q \varepsilon_p (q \times \overline{p}) \right] \quad \text{for } \mathbb{R}_1^3 \\ h(p,q) &= \frac{1}{2} \left[ -\varepsilon_p \varepsilon_q (p \times \overline{q}) - \varepsilon_q \varepsilon_p (q \times \overline{p}) \right] \quad \text{for } \mathbb{R}_2^4 \end{aligned} \qquad (2.3)$$

and it is called the semi-real quaternion inner product, [3]. The norm of a semi-real quaternion $q = (q_1, q_2, q_3, q_4) \in \mathbb{Q}_v$ is

$$N(q) = \sqrt{\left| q_1^2 + q_2^2 - q_3^2 - q_4^2 \right|}. \tag{2.4}$$

If $N(q) = 1$, then $q$ is called a semi-real unit quaternion, [5]. $q$ is called a semi-real spatial quaternion whenever $q + q = 0$, [2]. Moreover, quaternionic product of two semi-real spatial quaternions is $p \times q = g\ p, q\ + p \wedge_{\mathbb{L}} q$. $q$ is a semi-real temporal quaternion whenever $q - q = 0$. Any general $q$ can be written as $q = \frac{1}{2}(q + q) + \frac{1}{2}(q - q)$, [3].

## 3. Some Characterizations of Semi-Real Quaternionic Involute-Evolute Curves

The $4-$dimensional semi-Euclidean space $\mathbb{R}_2^4$ is identified with the space of unit semi-quaternions which is denoted by $\mathbb{Q}_v$. Let

$$\xi : I \subset \mathbb{R} \to \mathbb{Q}_v, \qquad s \to \xi(s) = \sum_{i=1}^{4} \xi_i(s) e_i, \qquad (1 \le i \le 4), \qquad e_4 = 1$$

be a smooth curve defined over the interval $I = [0,1]$. Let the arc-length parameter $s$ be chosen such that the tangent $\boldsymbol{T} = \xi'(s)$ has unit magnitude, [3]. Serret-Frenet apparatus of the semi-real quaternionic curve $\xi$ are given by

$$\begin{aligned}
\boldsymbol{T}(s) &= \xi', \\
\boldsymbol{N}(s) &= \frac{\xi''}{N(\xi'')}, \\
\boldsymbol{B}(s) &= \eta \varepsilon_n \varepsilon_T \left( \boldsymbol{E} \wedge_{\mathbb{L}} \boldsymbol{T} \wedge_{\mathbb{L}} \boldsymbol{N} \right), \\
\boldsymbol{E}(s) &= -\eta \varepsilon_n \varepsilon_b \varepsilon_T \varepsilon_N \frac{\boldsymbol{T} \wedge_{\mathbb{L}} \boldsymbol{N} \wedge_{\mathbb{L}} \xi'''}{N(\boldsymbol{T} \wedge_{\mathbb{L}} \boldsymbol{N} \wedge_{\mathbb{L}} \xi''')}, \qquad (\eta = \pm 1)
\end{aligned} \tag{3.1}$$

and

$$\begin{aligned}
\kappa(s) &= \varepsilon_N N(\xi''), \\
k(s) &= \varepsilon_N \frac{N(\boldsymbol{T} \wedge_{\mathbb{L}} \boldsymbol{N} \wedge_{\mathbb{L}} \xi''')}{N(\xi'')} \\
(r - \varepsilon_t \varepsilon_T \varepsilon_N \kappa)(s) &= \varepsilon_b \varepsilon_T \varepsilon_N \frac{h(\xi^{(iv)}, \boldsymbol{E})}{N(\boldsymbol{T} \wedge_{\mathbb{L}} \boldsymbol{N} \wedge_{\mathbb{L}} \xi''')}.
\end{aligned} \tag{3.2}$$

**Theorem 3.1.** Let $\{T, N, B, E\}$ be the Serret-Frenet frame in the point $\xi(s)$ of the semi-real quaternionic curve $\xi$ and $s$ is the arc-length parameter of the semi-real quaternionic curve $\xi$. Then the Serret-Frenet equations are

$$\begin{aligned}
T' &= \varepsilon_N \kappa N, \\
N' &= -\varepsilon_t \varepsilon_N \kappa T + \varepsilon_n k B, \\
B' &= -\varepsilon_t k N + \varepsilon_n \left(r - \varepsilon_t \varepsilon_T \varepsilon_N \kappa\right) E, \\
E' &= -\varepsilon_b \left(r - \varepsilon_t \varepsilon_T \varepsilon_N \kappa\right) B
\end{aligned} \tag{3.3}$$

where

$$\begin{aligned}
\kappa &= \varepsilon_N \|T'\|, \ N = \varepsilon_T (t \times T), \ B = \varepsilon_T (n \times T), \ E = \varepsilon_T (b \times T) \\
h(T,T) &= \varepsilon_T, \ h(N,N) = \varepsilon_N, \ h(B,B) = \varepsilon_n \varepsilon_T, \ h(E,E) = \varepsilon_b \varepsilon_T, \ [3].
\end{aligned} \tag{3.4}$$

It is obtained the Serret-Frenet formulae of the semi-real quaternionic curve $\xi = \xi(s)$ by making use of the Serret-Frenet formulae of the semi-real spatial quaternionic curve $\alpha = \alpha(s)$ where is $\alpha$ is a semi-real spatial quaternionic curve associated with the semi-real quaternionic curve $\xi$ and $\{t, n, b\}$ is the Frenet frame of the semi-real spatial quaternionic curve $\alpha$ in $\mathbb{R}_1^3$. Further, there are relationships between curvatures of the curves $\xi$ and $\alpha$. These relations are as follows: The torsion of $\xi$ is the principal curvature of the curve $\alpha$, the third curvature of $\xi$ is $(r - \varepsilon_t \varepsilon_T \varepsilon_N \kappa)$, where $r$ is the torsion of $\alpha$ and $\kappa$ is the principal curvature of $\xi$, [3].

**Definition 3.1.** Let $\phi, \xi : I \subset \mathbb{R} \to \mathbb{Q}_v$ be semi-real quaternionic curves with parameter $s^*$ and $s$, respectively. Moreover, $\{T_\phi, N_\phi, B_\phi, E_\phi\}$ and $\{T_\xi, N_\xi, B_\xi, E_\xi\}$ denote the Serret-Frenet frame of the curves $\phi$ and $\xi$, respectively. If

$$h(T_\phi(s^*), T_\xi(s)) = 0.$$

then, we call curves $\{\phi, \xi\}$ as semi-real quaternionic involute-evolute curves in $\mathbb{Q}_v$.

**Theorem 3.2.** Let $\xi, \phi : I \to \mathbb{Q}_v$ be a unit speed semi-real quaternionic curves. If semi-real quaternionic curve $\phi : I \to \mathbb{Q}_v$ is a quaternionic involute of the curve $\xi$, then we have $d_\mathbb{L}(\xi(s), \phi(s^*)) = |c - s|$, where $c$ is real number.

**Proof.** Let $\xi, \phi : I \to \mathbb{Q}_v$ be a unit speed semi-real quaternionic involute-evolute curves. From definition of semi-real quaternionic involute-evolute curves, we know that

$$\phi(s^*) = \xi(s) + \lambda(s)\mathbf{T}_\xi(s). \tag{3.5}$$

Then differentiating the equation (3.5) with respect to $s$, we get

$$\frac{d\phi}{ds^*}\frac{ds^*}{ds} = (1+\lambda')\mathbf{T}_\xi + \varepsilon_{N_\xi}\lambda\kappa_\xi \mathbf{N}_\xi.$$

Considering the last equation and using semi-real quaternion inner product with tangent vector $\mathbf{T}_\xi(s)$, we have

$$h\left(\frac{d\phi}{ds^*}\frac{ds^*}{ds}, \mathbf{T}_\xi(s)\right) = h\left(\left((1+\lambda')\mathbf{T}_\xi + \varepsilon_{N_\xi}\lambda\kappa_\xi \mathbf{N}_\xi\right), \mathbf{T}_\xi\right)$$

$$h\left(\frac{d\phi}{ds^*}\frac{ds^*}{ds}, \mathbf{T}_\xi(s)\right) = (1+\lambda')h(\mathbf{T}_\xi, \mathbf{T}_\xi) + \varepsilon_{N_\xi}\lambda\kappa_\xi h(\mathbf{N}_\xi, \mathbf{T}_\xi). \tag{3.6}$$

Considering the equation (3.4). we obtain that

$$\lambda(s) = c - s. \tag{3.7}$$

From the equations (3.5) and (3.7), we can reach

$$\phi(s^*) = \xi(s) + (c-s)\mathbf{T}_\xi(s). \tag{3.8}$$

Considering

$$N^2\left((c-s)\mathbf{T}_\xi\right) = \left|h\left((c-s)\mathbf{T}_\xi, (c-s)\mathbf{T}_\xi\right)\right| = (c-s)^2 \left|h(\mathbf{T}_\xi, \mathbf{T}_\xi)\right| = (c-s)^2 \left|\varepsilon_{T_\xi}\right| = (c-s)^2,$$

we obtain that

$$d_\mathbb{L}\left(\xi(s), \phi(s^*)\right) = N\left((c-s)\mathbf{T}_\xi\right) = |c-s|.$$

The following theorems provide some characterizations of semi-real quaternionic involute-evolute curves.

**Theorem 3.3.** Let $\phi, \xi : I \to \mathbb{Q}_v$ be unit speed semi-real quaternionic involute-evolute curves. The Serret-Frenet frame of semi-real quaternionic curve $\phi$, $\{\mathbf{T}_\phi, \mathbf{N}_\phi, \mathbf{B}_\phi, \mathbf{E}_\phi\}$ can be formed by frame of $\xi$, $\{\mathbf{T}_\xi, \mathbf{N}_\xi, \mathbf{B}_\xi, \mathbf{E}_\xi\}$.

**Proof.** Let us $\xi$ is a unit speed semi-real quaternionic curve. Without loss of generality, suppose that $\phi$ is the involute of $\xi$. By using the equation (3.3), if we calculate the derivatives of (3.8) with respect to $s$, we get

$$\mathbf{T}_\phi(s^*)\frac{ds^*}{ds} = \varepsilon_{N_\xi}(c-s)\kappa_\xi \mathbf{N}_\xi(s).$$

So, it is seen that

$$T_\phi(s^*) = \varepsilon_{N_\xi} N_\xi(s) \tag{3.9}$$

where $\dfrac{ds^*}{ds} = (c-s)\kappa_\xi$.

Let us calculate $N_\phi(s^*)$. If we calculate the derivatives of (3.9) with respect to $s$, we form that

$$T_\phi' = \phi'' = \varepsilon_{N_\xi} N_\xi'$$

$$\phi'' = -\varepsilon_t \kappa_\xi T_\xi + \varepsilon_n \varepsilon_{N_\xi} k B_\xi. \tag{3.10}$$

Thus by using the equation (3.4) and (3.10), and considering

$$N^2(\phi'') = |h(\phi'', \phi'')| = |k^2 h(B_\xi, B_\xi) + \kappa_\xi^2 h(T_\xi, T_\xi)| = |\varepsilon_n \varepsilon_{T_\xi} k^2 + \varepsilon_{T_\xi} \kappa_\xi^2|,$$

we obtain that

$$N(\phi'') = \sqrt{|\varepsilon_n \varepsilon_{T_\xi} k^2 + \varepsilon_{T_\xi} \kappa_\xi^2|}. \tag{3.11}$$

Using the equations (3.1), (3.10) and (3.11), we can write

$$N_\phi(s^*) = \frac{-\varepsilon_t \kappa_\xi T_\xi + \varepsilon_n \varepsilon_{N_\xi} k B_\xi}{\sqrt{|\varepsilon_n \varepsilon_{T_\xi} k^2 + \varepsilon_{T_\xi} \kappa_\xi^2|}} \tag{3.12}$$

Now let us calculate the binormal vector $E_\phi(s^*)$. Firstly, if we calculate derivate of the equation (3.10), we get

$$\phi''' = -\varepsilon_t \kappa_\xi' T_\xi - \varepsilon_t \varepsilon_{N_\xi}\left(\varepsilon_n k^2 + \kappa_\xi^2\right) N_\xi + \varepsilon_n \varepsilon_{N_\xi} k' B_\xi - \varepsilon_{N_\xi} k\left(r - \varepsilon_t \varepsilon_{T_\xi} \varepsilon_{N_\xi} \kappa_\xi\right) E_\xi. \tag{3.13}$$

Therefore, from the equations (3.9), (3.12) and (3.13) we can reach

$$T_\phi \wedge N_\phi \wedge \phi''' = \frac{-\varepsilon_n \varepsilon_{T_\xi} k\left(r - \varepsilon_t \varepsilon_{T_\xi} \varepsilon_{N_\xi} \kappa_\xi\right)\left(\varepsilon_{N_\xi} k T_\xi + \varepsilon_t \kappa_\xi B_\xi\right) + \varepsilon_t \varepsilon_n \varepsilon_b \varepsilon_{T_\xi}\left(\kappa_\xi' k - \kappa_\xi k'\right) E_\xi}{\sqrt{|\varepsilon_n \varepsilon_{T_\xi} k^2 + \varepsilon_{T_\xi} \kappa_\xi^2|}}.$$

(3.14)

The semi-real quaternionic norm of the equation (3.14) is

$$N(T_\phi \wedge N_\phi \wedge \phi''') = \sqrt{\left|\frac{\varepsilon_{T_\xi} k^2\left(r - \varepsilon_t \varepsilon_{T_\xi} \varepsilon_{N_\xi} \kappa_\xi\right)^2\left(k^2 + \varepsilon_n \kappa_\xi^2\right) + \varepsilon_b \varepsilon_{T_\xi}\left(\kappa_\xi' k - \kappa_\xi k'\right)^2}{\varepsilon_n \varepsilon_{T_\xi} k^2 + \varepsilon_{T_\xi} \kappa_\xi^2}\right|}.$$

(3.15)

Moreover, using the equations (3.14), (3.15) and considering $\varepsilon_n = \varepsilon_{n^*}\varepsilon_{b^*}$, $\varepsilon_{N_\phi} = \varepsilon_n \varepsilon_{N_\xi}$, we obtain that

$$\boldsymbol{E}_\phi = \eta \frac{\varepsilon_n \varepsilon_{T_\xi} k \left( r - \varepsilon_t \varepsilon_{T_\xi} \varepsilon_{N_\xi} \kappa_\xi \right) \left( \varepsilon_{N_\xi} k \boldsymbol{T}_\xi + \varepsilon_t \kappa_\xi \boldsymbol{B}_\xi \right) - \varepsilon_t \varepsilon_n \varepsilon_b \varepsilon_{T_\xi} \left( \kappa_\xi' k - \kappa_\xi k' \right) \boldsymbol{E}_\xi}{\sqrt{\left|\varepsilon_{T_\xi} k^2 \left( r - \varepsilon_t \varepsilon_{T_\xi} \varepsilon_{N_\xi} \right)^2 \left( k^2 + \varepsilon_n \kappa_\xi^2 \right) + \varepsilon_b \varepsilon_{T_\xi} \left( \kappa_\xi' k - \kappa_\xi k' \right)^2 \right|}}. \quad (3.16)$$

where $\eta = \pm 1$ providing that $\det \left( \boldsymbol{T}_\xi(s), \boldsymbol{N}_\xi(s), \boldsymbol{B}_\xi(s), \boldsymbol{E}_\xi(s) \right) = +1$. Similarly, considering the equations (3.1), (3.9), (3.12) and (3.16), if we make necessary arrangements, the binormal vector $\boldsymbol{B}_\phi(s^*)$ is obtained as follows

$$\boldsymbol{B}_\phi = \eta \varepsilon_{b^*} \varepsilon_n \varepsilon_b \varepsilon_{N_\xi} \frac{\left( \kappa_\xi k' - \kappa_\xi' k \right) \left( -\varepsilon_t k \boldsymbol{T}_\xi + \varepsilon_b \varepsilon_{N_\xi} \kappa_\xi \boldsymbol{B}_\xi \right) - \varepsilon_b \varepsilon_{N_\xi} k \left( r - \varepsilon_t \varepsilon_{T_\xi} \varepsilon_{N_\xi} \kappa_\xi \right) \left( k^2 + \varepsilon_n \kappa_\xi^2 \right) \boldsymbol{E}_\xi}{\sqrt{\left| \varepsilon_n \varepsilon_{T_\xi} k^2 + \varepsilon_{T_\xi} \kappa_\xi^2 \right|} \sqrt{\left| \varepsilon_{T_\xi} k^2 \left( r - \varepsilon_t \varepsilon_{T_\xi} \varepsilon_{N_\xi} \kappa_\xi \right)^2 \left( k^2 + \varepsilon_n \kappa_\xi^2 \right) + \varepsilon_b \varepsilon_{T_\xi} \left( \kappa_\xi' k - \kappa_\xi k' \right)^2 \right|}}.$$

(3.17)

So, the proof is completed.

**Theorem 3.4.** Let $\phi, \xi : I \to \mathbb{Q}_v$ be unit speed semi-real quaternionic involute-evolute curves. The Serret-Frenet curvatures of semi-real quaternionic curve $\phi$, $\left\{ \kappa_\phi, k^*, \left( r^* - \varepsilon_{t^*} \varepsilon_{T_\phi} \varepsilon_{N_\phi} \kappa_\phi \right) \right\}$ can be formed by curvatures of $\xi$, $\left\{ \kappa_\xi, k, \left( r - \varepsilon_t \varepsilon_{T_\xi} \varepsilon_{N_\xi} \kappa_\xi \right) \right\}$.

**Proof:** Without loss of generality, let $\phi, \xi : I \to \mathbb{Q}_v$ be unit speed semi-real quaternionic involute-evolute curves. Using the equations (3.2), (3.11) and considering $\varepsilon_{N_\phi} = \varepsilon_n \varepsilon_{N_\xi}$, we obtain that

$$\kappa_\phi = \varepsilon_n \varepsilon_{N_\xi} \sqrt{\left| \varepsilon_n \varepsilon_{T_\xi} k^2 + \varepsilon_{T_\xi} \kappa_\xi^2 \right|}. \quad (3.18)$$

Similarly, using the equations (3.4), (3.11) and (3.15), we have

$$k^* = \varepsilon_n \varepsilon_{N_\xi} \frac{\sqrt{\left| \varepsilon_{T_\xi} k^2 \left( r - \varepsilon_t \varepsilon_{T_\xi} \varepsilon_{N_\xi} \kappa_\xi \right)^2 \left( k^2 + \varepsilon_n \kappa_\xi^2 \right) + \varepsilon_b \varepsilon_{T_\xi} \left( \kappa_\xi' k - \kappa_\xi k' \right)^2 \right|}}{\left| \varepsilon_n \varepsilon_{T_\xi} k^2 + \varepsilon_{T_\xi} \kappa_\xi^2 \right|}. \quad (3.19)$$

Lastly, let us calculate the third curvature, $\left( r^* - \varepsilon_{t^*} \varepsilon_{T_\phi} \varepsilon_{N_\phi} \kappa_\phi \right)$. For this purpose if we calculate $\phi^{(iv)}$, we find that

$$\phi^{(iv)} = \left[-\varepsilon_t \kappa_\xi'' + \kappa_\xi \left(\varepsilon_n k^2 + \kappa_\xi^2\right)\right] T_\xi + \left[-3\varepsilon_t \varepsilon_{N_\xi} \kappa_\xi \kappa_\xi' - 3\varepsilon_t \varepsilon_n \varepsilon_{N_\xi} k k'\right] N_\xi$$
$$+ \left[-\varepsilon_t \varepsilon_n \varepsilon_{N_\xi} k \left(\varepsilon_n k^2 + \kappa_\xi^2\right) + \varepsilon_n \varepsilon_{N_\xi} k'' - \varepsilon_b \varepsilon_{N_\xi} k \left(r - \varepsilon_t \varepsilon_{T_\xi} \varepsilon_{N_\xi} \kappa_\xi\right)^2\right] B_\xi \quad (3.20)$$
$$+ \varepsilon_{N_\xi} \left[2k'\left(r - \varepsilon_t \varepsilon_{T_\xi} \varepsilon_{N_\xi} \kappa_\xi\right) + k\left(r - \varepsilon_t \varepsilon_{T_\xi} \varepsilon_{N_\xi} \kappa_\xi\right)'\right] E_\xi.$$

Considering the equation (3.20) and $E_\phi$ we get

$$h\left(\phi^{(iv)}, E_\phi\right) = \varepsilon_t \varepsilon_{N_\xi} \frac{\varepsilon_n k \left(r - \varepsilon_t \varepsilon_{T_\xi} \varepsilon_{N_\xi} \kappa_\xi\right)\left(\kappa_\xi'' k - \kappa_\xi k''\right) - \varepsilon_b \kappa_\xi k^2 \left(r - \varepsilon_t \varepsilon_{T_\xi} \varepsilon_{N_\xi} \kappa_\xi\right)^3}{\left[\left|\varepsilon_{T_\xi} k^2 \left(r - \varepsilon_t \varepsilon_{T_\xi} \varepsilon_{N_\xi} \kappa_\xi\right)^2 \left(k^2 + \varepsilon_n \kappa_\xi^2\right) + \varepsilon_b \varepsilon_{T_\xi} \left(\kappa_\xi' k - \kappa_\xi k'\right)^2\right|\right]}. \quad (3.21)$$

Lastly, using the equations (3.2), (3.15) and (3.21), we get the third curvature of the semi-real quaternionic curve $\phi$ as follows

$$r^* - \varepsilon_{t^*} \varepsilon_{T_\phi} \varepsilon_{N_\phi} \kappa_\phi = \varepsilon_{b^*} \varepsilon_t \varepsilon_n \varepsilon_{N_\xi} \frac{\varepsilon_n k\left(r-\varepsilon_t\varepsilon_{T_\xi}\varepsilon_{N_\xi}\kappa_\xi\right)\left(\kappa_\xi''k-\kappa_\xi k''\right)-\varepsilon_b\kappa_\xi k^2\left(r-\varepsilon_t\varepsilon_{T_\xi}\varepsilon_{N_\xi}\kappa_\xi\right)^3 + \left(2k'\left(r-\varepsilon_t\varepsilon_{T_\xi}\varepsilon_{N_\xi}\kappa_\xi\right)+k\left(r-\varepsilon_t\varepsilon_{T_\xi}\varepsilon_{N_\xi}\kappa_\xi\right)'\right)\left(\kappa_\xi'k-\kappa_\xi k'\right)}{\frac{1}{\sqrt{\left|\varepsilon_n\varepsilon_{T_\xi}k^2 + \varepsilon_{T_\xi}\kappa_\xi^2\right|}} \left[\left|\varepsilon_{T_\xi}k^2\left(r-\varepsilon_t\varepsilon_{T_\xi}\varepsilon_{N_\xi}\kappa_\xi\right)^2\left(k^2+\varepsilon_n\kappa_\xi^2\right)+\varepsilon_b\varepsilon_{T_\xi}\left(\kappa_\xi'k-\kappa_\xi k'\right)^2\right|\right]}.$$
$$(3.22)$$

This completed the proof.

From Theorem 3.3, we give following corollary.

**Corollary 3.1.** Let $\phi, \xi : I \to \mathbb{Q}_v$ be a unit speed semi-real quaternionic involute-evolute curves. If semi-real quaternionic curve $\xi$ is a semi-real quaternionic $w-curve$, then the Serret-Frenet frame of $\xi$ are

$$T_\phi(s^*) = \varepsilon_{N_\xi} N_\xi, \qquad N_\phi(s^*) = \frac{-\varepsilon_t \kappa_\xi T_\xi + \varepsilon_n \varepsilon_{N_\xi} k B_\xi}{\sqrt{\left|\varepsilon_n \varepsilon_{T_\xi} k^2 + \varepsilon_{T_\xi} \kappa_\xi^2\right|}}$$
$$B_\phi(s^*) = \eta \varepsilon_n \varepsilon_b \varepsilon_{b^*} E_\xi, \qquad E_\phi(s^*) = -\eta \varepsilon_n \varepsilon_{T_\xi} \frac{\varepsilon_{N_\xi} k T_\xi + \varepsilon_t \kappa_\xi B_\xi}{\sqrt{\left|\varepsilon_{T_\xi} k^2 + \varepsilon_n \varepsilon_{T_\xi} \kappa_\xi^2\right|}}.$$

(3.23)

**Theorem 3.5.** Let $\phi, \xi : I \subset \mathbb{R} \to \mathbb{Q}_v$ be semi-real quaternionic curves with parameter $s^*$ and $s$, respectively. If $(\phi, \xi)$ are semi-real quaternionic involute-evolute curves, then the semi-real spatial quaternionic curves $(\beta, \alpha)$, associated with $\phi$ and $\xi$, respectively, aren't the semi-real spatial quaternionic involute-evolute curves.

**Proof:** Let $(\phi, \xi)$ be semi-real quaternionic involute-evolute curves with parameter $s^*$ and $s$, respectively. The Frenet apparatus of $(\beta, \alpha)$, associated with $\phi$ and $\xi$, are $\{t^*, n^*, b^*, k^*, r^*\}$ and $\{t, n, b, k, r\}$, respectively. So, from the equation (3.4) we write $t^* = \varepsilon_{T_\phi} N_\phi \times T_\phi$ for the semi-real quaternionic curve $\phi = \phi(s^*)$. Here, if we use the equations (3,9) and (3.12), we can write $t^* = (xT_\xi + yB_\xi) \times \varepsilon_{N_\xi} N_\xi$, where

$$x = \frac{-\varepsilon_t \kappa_\xi}{\sqrt{\left|\varepsilon_{T_\xi} \kappa_\xi^2 + \varepsilon_n \varepsilon_{T_\xi} k^2\right|}}, \quad y = \frac{\varepsilon_n k}{\sqrt{\left|\varepsilon_{T_\xi} \kappa_\xi^2 + \varepsilon_n \varepsilon_{T_\xi} k^2\right|}}.$$

So, we obtain that $t^* = \varepsilon_{N_\xi} x(T_\xi \times N_\xi) + \varepsilon_{N_\xi} y(B_\xi \times N_\xi)$. By using the equations in (3.4), we find that

$$t^* = \varepsilon_{N_\xi} x\left(T_\xi \times \varepsilon_{T_\xi} T_\xi \times \hat{t}\right) + \varepsilon_{N_\xi} y\left(B_\xi \times \varepsilon_{T_\xi} T_\xi \times \hat{t}\right).$$

$$t^* = \varepsilon_{T_\xi} \varepsilon_{N_\xi} x\hat{t} + \varepsilon_{T_\xi} \varepsilon_{N_\xi} y\left(\varepsilon_{T_\xi} n \times \hat{t}\right)$$

$$t^* = -\varepsilon_{T_\xi} \varepsilon_{N_\xi} xt + \varepsilon_{N_\xi} yb.$$

From the last equation, we see that $t^*$ is perpendicular with $n$. However, if $t^*$ is linear depended with $n$, then $(\beta, \alpha)$ are semi-real spatial quaternionic involute-evolute curves. So, since $t^*$ is perpendicular with $n$, $(\beta, \alpha)$ aren't semi–real spatial quaternionic involute-evolute curves.

Now, we will give an example for the above theorem.

**Example 3.1.** We consider a quaternionic curve with the arc-length parameter $s$, $\xi : I \subset R \to R_2^4$ as noted below

$$\xi(s) = \left(\cosh(s), \sqrt{2}s, \sinh(s), \sqrt{2}\right)$$

for all $s \in I$. By considering the equations (3.1) and (3.2) we find that the Frenet apparatus of the quaternionic curve $\xi = \xi(s)$ as follows

$$T_\xi(s) = \left(\sinh(s), \sqrt{2}, \cosh(s), 0\right), \qquad N_\xi(s) = \left(\cosh(s), 0, \sinh(s), 0\right)$$

$$B_\xi(s) = -\varepsilon_b\left(\sqrt{2}\sinh(s), 1, \sqrt{2}\cosh(s), 0\right), \qquad E_\xi(s) = -\varepsilon_n \varepsilon_b(0, 0, 0, 1).$$

The curvatures functions of the quaternionic curve $\xi = \xi(s)$ are as follows

$$\kappa_\xi = -1, \quad k = \sqrt{2}, \quad r - \varepsilon_t \varepsilon_{T_\xi} \varepsilon_{N_\xi} \kappa_\xi = 0.$$

By using the equation (3.5), if we make necessary arrangement, we can easily find the quaternionic involute curve of the quaternionic curve $\xi = \xi(s)$, as follows

$$\phi(s) = \left((c-s)\sinh(s) + \cosh(s), \sqrt{2}c, (c-s)\cosh(s) + \sinh(s), \sqrt{2}\right)$$

which $c$ is a real number.

By using the equation (3.1), the spatial quaternionic curve $\alpha = \alpha(s)$ in $R_1^3$ associated with quaternionic curve $\xi = \xi(s)$ in $R_2^4$ is given by

$$\alpha(s) = \left(\sqrt{2}\sinh(s), s, -\sqrt{2}\cosh(s)\right)$$

where $s$ is the arc-length parameter of $\alpha$ and its curvature functions are as follows

$$k = \sqrt{2}, \quad r = -1.$$

Similarly, the spatial quaternionic curve $\beta$ in $R_1^3$ associated with quaternionic curve $\phi = \phi(s)$ in $R_2^4$ is given by $\beta(s) = (c, -s+c, c)$ which $c$ is a real number.

Now, if we calculate the quaternionic inner product $h(t \times t^*)$, then we obtain that $h(t \times t^*) = -1 \neq 0$. So, we can easily see that the spatial quaternionic curves $(\beta, \alpha)$ aren't spatial quaternionic involute-evolute curves.

The pictures of the some projections of the quaternionic curve $\xi = \xi(s)$, the quaternionic involute curves of $\xi$ and their associated spatial quaternionic curves are as follows

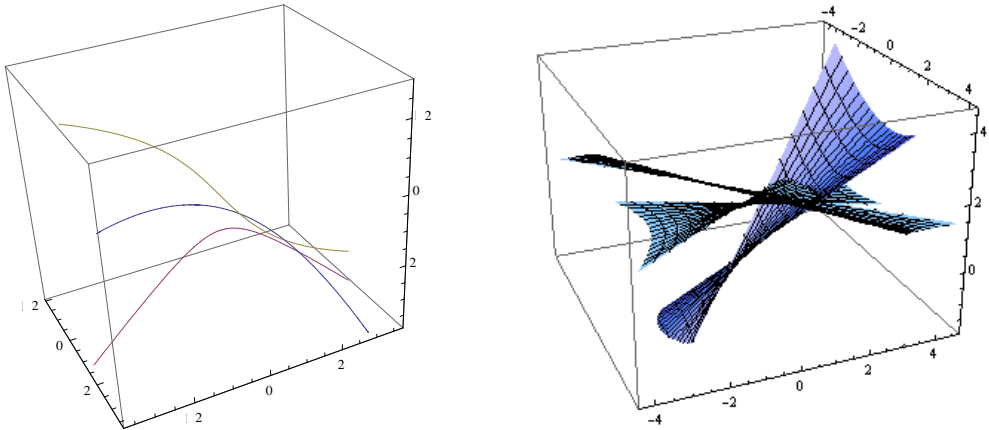

Figure 3.1. Some projections of the quaternionic curve $\xi = \xi(s)$ (on the left) and the quaternionic involutes of $\xi$ (on the right).

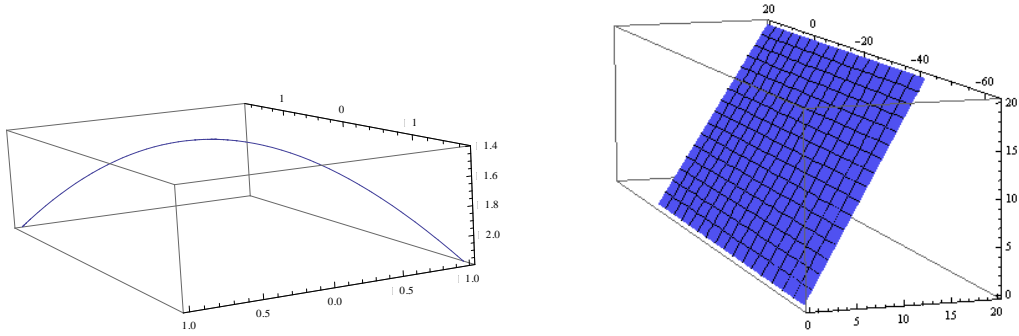

Figure 3.2. The spatial quaternionic curve associated with the quaternionic curve $\xi = \xi(s)$ (on the left) and the spatial quaternionic curves associated with the quaternionic involutes of $\xi$ (on the right).

*Department of Mathematics, Sakarya University, Sakarya, Turkey*
*E-mail: agungor@sakarya.edu.tr*

*Department of Mathematics, Sakarya University, Sakarya, Turkey*
*E-mail: tsoyfidan@sakarya.edu.tr*